\documentclass{amsart}
 
 \usepackage{amscd}
 \usepackage{amssymb}

\usepackage{amsthm}

\newtheorem{lem}[equation]{Lemma}

\newtheorem{prop}[equation]{Proposition}

\newtheorem*{thm*}{Theorem}
\newtheorem*{prop*}{Proposition}
\newtheorem*{cor*}{Corollary}
\newtheorem*{lem*}{Lemma}
\newtheorem*{MT*}{Main Theorem}

\newtheorem{specthm}{Theorem}

\theoremstyle{definition} %

\newtheorem*{defn*}{Definition}

\newtheorem{eg}[equation]{Example}

\theoremstyle{remark} %
\newtheorem{rmk}[equation]{Remark}

\newtheorem*{rmk*}{Remark}
\newtheorem*{rmks*}{Remarks}

\newtheoremstyle{exercise}
  {3pt}
  {3pt}
  {\small}
  {}
  {\sc\small}
  {.}
  {.5em}
   {}     
  {}

\theoremstyle{exercise}

\makeatletter
{\renewcommand{\theequation}{#1}}%
{\renewcommand{\theequation}{\arabic{equation}}\addtocounter{equation}{-1}\global\@ignoretrue}
\makeatother

\makeatletter
{\renewcommand{\theequation}{#1}\begin{eqnarray}}%
{\end{eqnarray}\renewcommand{\theequation}{\arabic{equation}}\addtocounter{equation}{-1}\global\@ignoretrue}
\makeatother

\makeatletter
{\smallskip \refstepcounter{equation}\noindent{\textbf{\theequation.} }{{\textbf{#1.}}}}%
{\smallskip \global\@ignoretrue}
\makeatother

\makeatletter
{\smallskip \refstepcounter{equation}{\sc \theequation}{\sc (#1).}}%
{\smallskip \global\@ignoretrue}
\makeatother

\makeatletter
{\smallskip \refstepcounter{equation}\noindent{\sc \theequation.}{\sl{ #1.}}}%
{\smallskip \global\@ignoretrue}
\makeatother

\makeatletter
\newenvironment{borel*}%
{\smallskip \refstepcounter{equation}\noindent{\textbf{\theequation.}}}%
{\global\@ignoretrue}
\makeatother

\newcommand{\flist}[1]{\hangindent\leftmargini\textup{(1)}\hskip\labelsep {#1}%
\begin{enumerate}%
\setcounter{enumi}{1}%
}

\newcommand{\ot}{\otimes}

\newcommand{\eand}{\quad\text{and}\quad}

\newcommand{\Q}{{\mathbb{Q}}}        
\newcommand{\F}{{\mathbb{F}}}        
\newcommand{\Z}{{\mathbb{Z}}}        
\renewcommand{\H}{{\mathcal{H}}}  

\newcommand{\la}{\lambda}

\newcommand{\oddots}{{\mathinner{\mkern1mu\raise1pt\vbox{\kern7pt\hbox{.}}\mkern2mu\raise4pt\hbox{.}\mkern2mu\raise7pt\hbox{.}\mkern1mu}}}

\newcommand{\ksep}{k_{{\mathrm{sep}}}}

\newcommand{\kx}{k^\times}

\newcommand{\qform}[1]{{\left\langle{#1}\right\rangle}}                   

\DeclareMathOperator{\SL}{SL}
\DeclareMathOperator{\GL}{GL}

\DeclareMathOperator{\Int}{Int}
\DeclareMathOperator{\chr}{char}

\DeclareMathOperator{\End}{End}

\newcommand{\ra}{\rightarrow}

\newcommand{\stcolvec}[2]{\left( \begin{smallmatrix} #1 \\ #2 \end{smallmatrix} \right) }
\newcommand{\stbtmat}[4]{\left( \begin{smallmatrix} #1&#2 \\ #3&#4 \end{smallmatrix} \right) }

 \theoremstyle{plain}
 \newtheorem*{specA}{Theorem A} 
 \newtheorem*{specB}{Theorem B} 
 
 \numberwithin{equation}{section}
 
 \newcommand{\sbar}{\bar{s}}

\DeclareMathOperator{\PGL}{PGL}

\newcommand{\even}{{{\mathrm{even}}}}
\newcommand{\odd}{{{\mathrm{odd}}}}

\renewcommand{\F}{\mathbb{F}}

\newcommand{\cprod}{{\mathinner{\mkern2mu\raise2.5pt\hbox{.}}}}

\newcommand{\rbar}{\overline{\rho}}

\newcommand{\nd}{{\mathrm{nd}}}

\newcommand{\podd}{\phi^{{\mathrm{odd}}}_n}
\newcommand{\peven}{\phi^{{\mathrm{even}}}_n}

\begin{document}
 \title{Orthogonal representations of twisted forms of $\SL_2$}
 
 \author{Skip Garibaldi}
\address{Department of Mathematics \& Computer Science, Emory University, Atlanta, GA 30322, USA}
\email{skip@member.ams.org}
\urladdr{http://www.mathcs.emory.edu/{\textasciitilde}skip/}


\thanks{{\tt Version of \today.}}
\thanks{The author's research was partially supported by the NSF}

\begin{abstract}
For every absolutely irreducible orthogonal representation of a twisted form of $\SL_2$ over a field of characteristic zero, we compute the ``unique" symmetric bilinear form that is invariant under the group action.  We also prove the analogous result for Weyl modules in prime characteristic (including characteristic 2) and an isomorphism between two symmetric bilinear forms given by binomial coefficients.
\end{abstract}

\maketitle

\setlength{\unitlength}{.75cm}

Let $V$ be an absolutely irreducible and finite-dimensional representation of a group $G$ over a field $k$.  If there is a nonzero quadratic form $q$ on $V$ that is invariant under $G$, then by Schur's Lemma $q$ is uniquely determined up to multiplication by an element of $\kx$.  It is natural to ask: \emph{What is $q$?} 
On the one hand, the form $q$ gives some information about $G$.  For example, when $V$ is the adjoint representation, $q$ is essentially the Killing form on the Lie algebra of $G$.  On the other hand, if $G$ is contained in a group $H$, then by restricting representations of $H$, one reduces the problem of computing the $H$-invariant quadratic forms to the problem of computing those invariant under $G$.  This is the method used in \cite{Jac:ex} to calculate the Killing forms of various Lie algebras.

The purpose of this paper is to answer the question ``What is $q$?" for the smallest groups, i.e., groups of type $A_1$, meaning twisted forms of $\SL_2$ and $\PGL_2$.
We study Weyl modules $V$ instead of absolutely irreducible representations.  (If $k$ has characteristic zero, this is no change at all.   In prime characteristic, the results for Weyl modules imply those for irreps, see \ref{VvL}.)
%

Recall that the group $\SL_2$ has a unique fundamental dominant weight $\omega$.  We write $V(n)$ for the Weyl module with highest weight $n\omega$, where $n$ is a positive integer.  If $k$ has characteristic zero, $V(n)$ is isomorphic to the $n$-th symmetric power of $k^2$.  For every field $k$, every group of type $A_1$, and $n$ even, the representation $V(n)$ is $k$-defined (as can easily be seen by Galois descent as in \cite[p.~139]{BoTi}) and has a nonzero quadratic form invariant under the group action.  We compute this quadratic form.  For remarks on the modules $V(n)$ with $n$ odd, see \ref{V.odd}.

\smallskip

We now state our main results in the case where $k$ has characteristic $\ne 2$.
The situation in characteristic 2 is very different from other characteristics, and \emph{from here through \S\ref{rep.pf}, we only consider fields of characteristic $\ne 2$.}  

Define a symmetric bilinear form $\podd$ by
\begin{equation} \label{odd.def}
\podd := \bigoplus_{\stackrel{0 \le i < n/2}{\text{$i$ odd}}} \qform{ \binom{n}{i} }.
\end{equation}
The notation $\qform{\la}$ denotes a one-dimensional vector space spanned by  a vector $e$ endowed with a symmetric bilinear form that sends $(e,e) \mapsto \la$ for $\la \in \kx$ and $\oplus$ means `orthogonal sum'.  The form $\phi^\odd_2$ is 0-dimensional, and $\phi^\odd_4$ is isomorphic to $\qform{1}$. We refer to \cite{Lam} for general background on symmetric bilinear forms and quaternion algebras.  In particular, since the characteristic is not 2, we can and do identify symmetric bilinear forms with quadratic forms by sending a symmetric bilinear form $b$ to the quadratic form $v \mapsto b(v,v)$.

Every semisimple algebraic group of type $A_1$ over $k$ is isogenous to $\SL(Q)$ for some quaternion $k$-algebra $Q$, where $\SL(Q)$ denotes the group whose $k$-points are the norm 1 elements of $Q$.  We write $Q$ also for the norm form on $Q$; it is a quadratic form.  We write $Q'$ for the unique 3-dimensional quadratic form such that $Q = \qform{1} \oplus Q'$.

\begin{specthm}[$\chr k \ne 2$] \label{REP}
For $n$ even, there is an $\SL(Q)$-invariant symmetric bilinear form on the Weyl module $V(n)$ that is isomorphic to 
\[
\qform{2  } \cdot \podd \cdot Q \oplus \begin{cases}
\qform{\binom{n}{n/2}} & \text{if $n \equiv 0 \bmod{4}$}\\
\qform{\binom{n}{n/2}}Q' & \text{if $n \equiv 2 \bmod{4}$.}
\end{cases}
\]
\end{specthm}

The smallest case of Theorem \ref{REP}\ is where $n = 2$, i.e., the adjoint representation of $\SL(Q)$.  There we find the symmetric bilinear form $\qform{2}Q'$, whereas the Killing form is $\qform{-2} Q'$.

\smallskip

The proof of Th.~\ref{REP}\ gives a result on symmetric bilinear forms as a byproduct.  Define a form $\peven$ as in \eqref{odd.def}, except with the word ``odd" replaced with ``even" throughout.  For example, $\phi^\even_2$ and $\phi^\even_4$ are $\qform{1}$, and $\phi^\even_{6}$ is $\qform{1, 15}$.  
\begin{specthm}[$\chr k \ne 2$] \label{QUAD}
We have the following isomorphism of symmetric bilinear forms over $k$:
\[
\peven \cong \begin{cases}
\podd & \text{if $n \equiv 0 \bmod{4}$} \\
\qform{2 \binom{n}{n/2}} \oplus \podd & \text{if $n \equiv 2 \bmod{4}$}.
\end{cases}
\]
\end{specthm}

As an illustration, for $n = 16$, Th.~\ref{QUAD}\ asserts the isomorphism
\[
\qform{1, 2^3 \cdot 3 \cdot 5, 2^2 \cdot 5 \cdot 7 \cdot 13, 2^3 \cdot 7 \cdot 11 \cdot 13} \cong \qform{2^4, 2^4 \cdot 5 \cdot 7, 2^4 \cdot 3 \cdot 7 \cdot 13, 2^4 \cdot 5 \cdot 11 \cdot 13}.
\]
  Note that the forms $\peven$ and $\podd$ may be degenerate if the characteristic of $k$ is not zero.

One imagines that Theorem \ref{QUAD}\ is known, but the experts I consulted said they had never seen it before.  For small $n$ and fields of characteristic zero, one can check the theorem using the Hasse-Minkowski Principle.  The proof of Theorem \ref{QUAD}\ given below uses neither Hasse-Minkowski nor clever manipulations of binomial coefficients.

\smallskip

The text through \S\ref{rep.pf} is concerned with proving Theorems \ref{REP} and \ref{QUAD} in characteristic different from 2.  Section \ref{quads} discusses orthogonal representations over all fields, with a special focus on the situation in characteristic 2.  The final section, \S\ref{char2}, states and proves analogues of Theorems \ref{REP} and \ref{QUAD} in characteristic 2.  We remark that the version of Th.~\ref{REP} in characteristic 2 is very different from the result in other characteristics, in that in characteristic 2 the only Weyl module whose invariant quadratic form depends on $Q$ is $V(2)$.

\section{Binomial coefficients} \label{binomial}

This section is purely about integers, and we write $n$ for a positive integer.  Kummer proved that the power of a prime $p$ dividing the binomial coefficient $\binom{n}{m}$ is the number of carries when you add $m$ and $n - m$ in base $p$.  We will use the following special case.
For a non-negative integer $m$, the base-$p$ expansion of $m$ is the list of integers $m_0, m_1, \ldots$ such that $0 \le m_j < p$ for all $j$ and $m = \sum_j m_j p^j$.
\begin{lem}[Kummer, 1852] \label{Kummer}
For $0 \le m \le n$, the binomial coefficient $\binom{n}{m}$ is divisible by a prime $p$ if and only if 
in the base-$p$ expansions of $m$ and $n$, there is some place $j$ with $m_j > n_j$. 
\end{lem}

\begin{proof} 
See, for example, \cite[Th.~1]{Fine} or \cite{Granville}.
 \end{proof}

\begin{lem} \label{binom.lem}
Let $p$ be an odd prime.  If $n$ is divisible by $4$, the lists
\[
\binom{n}{0}, \binom{n}{2}, \ldots, \binom{n}{n/2-2} \eand
\binom{n}{1}, \binom{n}{3}, \ldots, \binom{n}{n/2-1}
\]
have the same number of elements that are divisible by $p$.  

If $n$ is congruent to $2 \bmod{4}$, then the lists 
\[
\binom{n}{0}, \binom{n}{2}, \ldots, \binom{n}{n/2-1}
\eand
\binom{n}{1}, \binom{n}{3}, \ldots, \binom{n}{n/2-2}, \binom{n}{n/2}
\]
have the same number of elements that are divisible by $p$.
\end{lem}

The analogous statement is false for $p = 2$, see the characteristic 2 version of Th.~\ref{QUAD} in \S\ref{char2} below.

\begin{proof}
We write $n_j$ and $m_j$ for the $j$-th terms in the base-$p$ expansions of $n$ and $m$.  

\smallskip
{\emph{\underline{Case 1.}}} Suppose first that some $n_\ell$ of $n$ is odd.  
Whether $n$ is congruent to 0 or 2 mod 4, the two lists have the same length.  So it suffices to prove that the lists have the same number of elements that are \emph{not} divisible by $p$.

Write $h_0, h_1, \ldots$ for the base-$p$ expansion of $n/2$, so $n/2 = \sum h_j p^j$ and $n = \sum (2 h_j) p^j$.  But since $n_\ell$ is odd, $2h_0, 2h_1, \ldots$ cannot be the base $p$ representation of $n$, and $p$ divides $\binom{n}{n/2}$.

Let $S$ be the set of numbers whose base-$p$ expansion is $\sum s_j p^j$ with $0 \le s_j \le n_j$ for all $j$ and $s_\ell$ is zero.  Write $\hat{S}$ for the collection of numbers
\begin{equation}  \label{s.list}
s, s + p^\ell, s+2p^\ell, \ldots, s+n_\ell p^\ell
\end{equation}
as $s$ varies over elements of $S$.  The set $\hat{S}$ is precisely the set of numbers $m$ with $0 \le m \le n$ such that $p$ does not divide $\binom{n}{m}$.

We note that half of the numbers in \eqref{s.list} are even and half are odd, since $n_\ell$ and $p$ are both odd.  As $\binom{n}{m}$ equals $\binom{n}{n-m}$, the numbers $m, n-m$ have the same parity, and $n/2$ is not in $\hat{S}$, we conclude that one quarter of $\hat{S}$ belongs to each of the two lists. This proves the lemma in this case.

\smallskip
{\emph{\underline{Case 2.}}}
Now suppose that every base-$p$ digit $n_j$ of $n$ is even.
 We examine the set $\hat{T}$ of numbers $m$ such that $0 \le m \le n$ and $p$ divides $\binom{n}{m}$.  Such $m$'s are exactly those whose base-$p$ expansion $m_0, m_1, \ldots$ has a position $\ell$ such that $m_j = n_j$ for $j > \ell$,\footnote{If $\ell$ is such that $n_j = 0$ for $j > \ell$, then this condition just says that $m_j = 0$, which is a consequence of the hypothesis that $m \le n$.} $m_\ell < n_\ell$, and there is some $j < \ell$ such that $m_j > n_j$.  Put $T_\ell$ for the set of integers $t$ whose base-$p$ expansion $t_0, t_1, \ldots$ has $t_j = n_j$ for $j > \ell$, has $t_\ell = 0$, and has $t_j > n_j$ for some $j < \ell$.  Then $\hat{T}$ is the union of the (disjoint) collections
\begin{equation} \label{t.list}
t, t+p^\ell, t+2p^\ell, \ldots, t+(n_\ell - 1)p^\ell
\end{equation}
for $t \in T_\ell$ and $0 < \ell \le r$.
Half of \eqref{t.list} is even and half is odd, so half of $\hat{T}$ is even and half is odd.  Further $\hat{T}$ does not contain $\binom{n}{n/2}$ by the hypothesis on $n$.  As in Case 1, we conclude that the two lists have the same number of terms that are divisible by $p$.
\end{proof}

\section{Description of $V(n)$ in the split case} \label{split}

Suppose that $Q$ is split, i.e., that $\SL(Q)$ is the split group $\SL_2$.  Recall the well-known $\SL_2(\Q)$-invariant bilinear form on the $n$-th symmetric power $S^n(\Q^2)$:  Fix a basis $e_0, e_1, \ldots, e_n$ given by
\[
e_m := \binom{n}{m} \left[ \begin{smallmatrix} 1 \\ 0 \end{smallmatrix} \right]^{n-m}
\left[ \begin{smallmatrix} 0 \\ 1 \end{smallmatrix} \right]^{m}
\]
where $\binom{n}{m}$ denotes a binomial coefficient and $\left[ \begin{smallmatrix} 1 \\ 0 \end{smallmatrix} \right]$ and $\left[ \begin{smallmatrix} 0 \\ 1 \end{smallmatrix} \right]$ are the standard basis vectors for $\Q^2$.  The bilinear form $f$ on $S^n(\Q^2)$ given by
\begin{equation} \label{f.def}
f(e_\ell, e_m) = \begin{cases}
0 & \text{if $\ell + m \ne n$} \\
(-1)^\ell \binom{n}{\ell} & \text{if $\ell + m = n$}
\end{cases}
\end{equation}
is preserved by $\SL_2$  \cite[\S{VIII.1.3}, Remark 3]{Bou:g7}.

We view the group $\SL_2$ as a Chevalley group defined over $\Z$ as in \cite{St} by fixing a pinning.  The elements $e_0, e_1, \ldots, e_n$ in $S^n(\Q^2)$ span a minimal admissable lattice $L$ in $S^n(\Q^2)$, see e.g.~\cite[\S27, Exercise 8]{Hum:LA}.  The Weyl module $V(n)$ with highest weight $n\omega$ over the field $k$ is $L \ot_\Z k$.  In particular, it has basis $e_0, e_1, \ldots, e_n$, and the form $f$ extends to an $\SL_2$-invariant bilinear form on the $k$-vector space $V(n)$.

\begin{borel*} 
The representation $V(n)$ is irreducible if and only if 
\begin{equation} \label{ndeg}
\parbox{4in}{$k$ has characteristic zero or $k$ has prime characteristic $p$ and $n+1 = rp^s$ for integers $r, s$ with $0 \le r < p$,}
\end{equation}
see e.g.~\cite[p.~240]{WinterSL2}.  By Kummer (or \S\ref{quads} below), this condition is also equivalent to $f$ being nondegenerate over $k$, i.e., that none of the binomial coefficients $\binom{n}{m}$ are divisible by the characteristic for $1 \le m \le n$.  
\end{borel*}

From here on, we restrict to the case where \emph{$n$ is even}.

\begin{rmk}
Within the line of $\SL_2$-invariant symmetric bilinear forms on $V(n)$, the isomorphism class of $f$ is distinguished by the fact that
\[
f(e_0, \rho\stbtmat{0}{1}{-1}{0} e_0) = f(e_0, (-1)^n e_0) = 1
\]
is a square in $\kx$.  More generally, for every maximal torus $T$ contained in a Borel subgroup $B$ of $\SL_2$, every highest weight vector $v \in V(n)$ (relative to $T$ and $B$), and every element $w \in N_{\SL_2}(T)(k)$ that represents the element $-1$ in the Weyl group, we have: \emph{If $s$ is an $\SL_2$-invariant symmetric bilinear form on $V(n)$ such that $s(v, wv)$ is a square in $\kx$, then $s$ is isomorphic to $f$.}  This follows from the conjugacy of pairs $T \subset B$ and the fact that $n$ is even.
\end{rmk}

\section{Computations}\label{desc}

We now address Theorem \ref{REP}\ in the case where $Q$ is not split.  We maintain the assumptions that $n$ is even and the characteristic of $k$ is different from 2.

The representation $\rho \!: \SL_2 \ra \GL(V(n))$ kills the center $\pm 1$ of $\SL_2$ and so induces a homomorphism $\rbar \!: \PGL_2 \ra \GL(V(n))$, where $\PGL_2$ denotes the (algebraic-group-theoretic) quotient $\SL_2/\{ \pm 1 \}$.  As the form $f$ defined in \S\ref{split} is $\SL_2$-invariant, the images of $\rho$ and $\rbar$ lie in the orthogonal group $O(f)$ of $f$, and we obtain a map in Galois cohomology
\begin{equation} \label{Oinv}
H^1(k, \PGL_2) \ra H^1(k, O(f)).
\end{equation}
The sets $H^1(k, \PGL_2)$ and $H^1(k, O(f))$ classify $k$-forms of $\SL_2$ and symmetric bilinear forms over $k$ that become isomorphic to $f$ over a separable closure of $k$, respectively \cite[\S{III.1}]{SeCG}, and the image of a group $\SL(Q)$ in $H^1(k, O(f))$ is an $\SL(Q)$-invariant symmetric bilinear form $f_n$ on $V$.

\begin{rmk}  
The map \eqref{Oinv} is defined for every extension $K/k$ and so leads to a morphism of functors $H^1(*, \PGL_2) \ra H^1(*, O(f))$. Suppose that $f$ is nondegenerate, i.e., \eqref{ndeg} holds. By \cite[28.4]{Se:ci}, 
the map \eqref{Oinv} is necessarily of the form
\[
Q \mapsto \psi_0 \oplus \psi_1 \cdot Q'
\]
for uniquely-determined symmetric bilinear forms $\psi_0, \psi_1$ defined over $k$ with $\psi_1$ anisotropic.  The point of Th.~\ref{REP}\ is to give $\psi_0$ and $\psi_1$ explicitly.
\end{rmk}

The $k$-algebra $Q$ is generated by elements $i, j$ such that $i^2 = a$, $j^2 = b$, and $ij = -ji$ for some $a, b \in \kx$.  As $Q$ is not split, $a$ is not a square in $k$ and $k(\sqrt{a})$ is a proper extension.  The algebra $Q$ is isomorphic to the $k$-subalgebra of $M_2(k(\sqrt{a}))$ generated by
\[
\stbtmat{\sqrt{a}}{0}{0}{-\sqrt{a}} \quad \text{and} \quad \stbtmat{0}{b}{1}{0},
\]
which we may identify with $i$ and $j$ respectively.
For $\iota$ the nonidentity $k$-automorphism of $k(\sqrt{a})$, $Q$ is obtained from $M_2(k)$ by twisting by the cocycle $\eta \in Z^1(k(\sqrt{a})/k, \PGL_2)$ given by
\begin{equation} \label{eta.def}
\eta_\iota = \Int \stbtmat{0}{b}{1}{0},
\end{equation}
where we write $\Int(g)$ for the map $x \mapsto gxg^{-1}$ on $M_2$.  (Recall that  for an extension $K/k$, the $K$-points of $\PGL_2$ are of the form $\Int(g)$ for $g \in \GL_2(K)$.)

We now compute $\rbar(\eta_\iota) \in O(f)(k(\sqrt{a}))$.  The element
\[
\beta := \stbtmat{0}{-\sqrt{-b}}{1/\sqrt{-b}}{0} \quad \in \SL_2(\ksep)
\]
is a scalar multiple of $\eta_\iota$, hence has the same image in $\PGL_2(k)$ and $\rbar(\eta_\iota)$ equals $\rho(\beta)$.  We find:
\[
\rbar(\eta_\iota) e_m =  \rho(\beta) e_m = (-1)^{n/2}b^{m-(n/2)} e_{n-m}.
\]
The action of $\rbar(\eta_\iota)$ preserves the line spanned by $e_{n/2}$ and the hyperbolic planes spanned by $e_m, e_{n -m}$ for $m \ne n/2$.  We compute separately the restriction of $f_n$ to each of these subspaces.

For the line spanned by $e_{n/2}$, we have $\rbar(\eta_\iota) e_{n/2} = (-1)^{n/2} e_{n/2}$, so the fixed $k$-subspace of the line is spanned by $e_{n/2}$ if $n/2$ is even and by $\sqrt{a} e_{n/2}$ if $n/2$ is odd.  Therefore the restriction of $f_n$ to this subspace is 
\[\qform{ (-a)^{n/2} \binom{n}{n/2}}.\]

For the $k(\sqrt{a})$-plane spanned by $e_m, e_{n-m}$ with $m \ne n/2$, the $\rbar(\eta_\iota) \iota$-fixed $k$-subspace is spanned by
\[
e_m + (-1)^{n/2}b^{m -(n/2)} e_{n-m} \quad \text{and} \quad \sqrt{a} e_m - (-1)^{n/2} b^{m-(n/2)}  \sqrt{a} e_{n-m}.
\]
These vectors are orthogonal and the restriction of $f_n$ to this subspace is 
\[
\qform{2} \qform{1, -a} \qform{(-b)^{m-(n/2)} \binom{n}{m}}.
\]

Summarizing the previous two paragraphs, we find that $f_n$ is isomorphic to
\[
\qform{2} \qform{1, -a} \left[ \bigoplus_{\stackrel{0 \le m < n/2}{\text{$m - (n/2)$ odd}}} \qform{-b} \qform{\textstyle\binom{n}{m}} \oplus \bigoplus_{\stackrel{0 \le m < n/2}{\text{$m - (n/2)$ even}}} \qform{\textstyle\binom{n}{m}} \right] \oplus \qform{(-a)^{n/2} \textstyle\binom{n}{n/2}}.
\]
That is:
\begin{equation} \label{desc.summ}
\parbox{4in}{\emph{If $n$ is divisible by $4$, then $f_n$ is isomorphic to 
\[
\qform{2} \qform{1,-a} \left[ \qform{-b} \podd \oplus \peven \right] \oplus \qform{\binom{n}{n/2}}.
\]
If $n$ is congruent to $2$ mod $4$, then $f_n$ is isomorphic to
\[
\qform{2} \qform{1, -a} \left[ \qform{-b} \peven \oplus \podd \right] \oplus \qform{-a\binom{n}{n/2}}.
\]
}}
\end{equation}

\section{Proof of Theorem \ref{QUAD}}

This purpose of this section is to prove Theorem B.  We put $q := \podd$ if 4 divides $n$ and $q := \podd \oplus \qform{2 \binom{n}{n/2}}$ if $n$ is congruent to 2 mod 4; we want to show that $\peven$ and $q$ are isomorphic.

We consider the ``generic" case of Th.~\ref{REP}, i.e., where $k$ is $\F(a, b)$ where $\F$ is a prime field (of characteristic $\ne 2$) and $a, b$ are independent indeterminates.  
We write $f_n$ for the form computed in \eqref{desc.summ}.  Recall that a symmetric bilinear form $s$ on a vector space $V$ induces a nondegenerate form --- the \emph{nondegenerate part} of $s$
--- on $V$ modulo the radical of $s$.
By Springer's Theorem (twice), the coefficients of $\qform{a}$ and $\qform{b}$ in the nondegenerate part of $f_n$ --- namely, the nondegenerate parts of $\qform{-2} \peven$ and $\qform{-2} q$ or vice versa --- are determined up to Witt equivalence.
The isomorphism class of the form $f_n$ depends only on the isomorphism class of $Q$ and not on the choice of $\eta$.  In particular, $f_n$ is unchanged when we interchange the roles of $a$ and $b$, and we conclude that $\peven$ and $q$ have Witt-equivalent nondegenerate parts.  By Lemma \ref{binom.lem},
the forms $\peven$ and $q$ have nondegenerate parts of the same dimension.  Since the total dimension of $\peven$ and $q$ is also the same, we conclude that the two forms are isomorphic.  This completes the proof of Th.~\ref{QUAD}.$\hfill\qed$

\section{Proof of Theorem \ref{REP}} \label{rep.pf}

We now prove Theorem \ref{REP} for fields of characteristic $\ne 2$.  In the case where $Q$ is not split, one combines Theorem \ref{QUAD}\ with \eqref{desc.summ} and the fact that the norm of $Q$ is $\qform{1, -a, -b, ab}$.  Therefore, we assume that $Q$ is split, in which case the form $f$ from \eqref{f.def} is an invariant quadratic form on $V(n)$.
For $0 \le m < n/2$, the elements $e_m, e_{n-m}$ in $V(n)$ span a hyperbolic plane $\H$ if $\binom{n}{m}$ is nonzero and a 2-dimensional subspace of the radical otherwise.  
We conclude that $f$ is isomorphic to
\[
\left[ \H \cdot (\peven \oplus \podd) \right] \oplus \qform{(-1)^{n/2} \binom{n}{n/2}}.
\]
(Specialists in symmetric bilinear forms are used to ignoring hyperbolic planes.  Here we must keep careful track, because some are multiplied by zero.) 
By Theorem \ref{QUAD}, this is isomorphic to 
\[
\left[ (\H \oplus \H) \cdot \podd \right]  \oplus 
\begin{cases}
\qform{ \binom{n}{n/2}} & \text{if $n \equiv 0 \bmod{4}$} \\
\qform{\binom{n}{n/2}} \left[ \qform{- 1} \oplus \H  \right]& \text{if $n \equiv 2 \bmod{4}$.}
\end{cases}
\]
As $Q$ is split, its quadratic norm form is $\H \oplus \H$ and $Q'$ is isomorphic to $\qform{-1, 1, -1}$.  This completes the proof of Th.~\ref{REP}.$\hfill\qed$

\section{Invariant quadratic forms on Weyl modules} \label{quads}

In this section, $k$ is any field (including the possibility of characteristic 2) and $G$ is a split semisimple linear algebraic group over $k$.  We extend some results on $G$-invariant bilinear forms on an irreducible representation $L(\la)$ with highest weight $\la$ to the case of a Weyl module $V(\la)$.  Although the extensions are straightforward, we write them explicitly for lack of a better reference.  The results are mostly interesting in characteristic 2, but we state them as much as possible in a characteristic-free way.

\begin{borel*} \label{clar}
First, a clarification.  We say that a form $f$ on a representation $V$ of $G$ over  $k$ is \emph{$G$-invariant} if $f$ is fixed by the algebraic group $G$ (and not merely by the group of $k$-points $G(k)$).  As $G$ is reductive, this is equivalent to saying that $G(K)$ fixes the form $f$ on $V \ot K$ for some infinite field $K$ containing $k$.  As an example of the difference,  the form $\stcolvec{x}{y} \mapsto x^2 + xy + y^2$ on $(\mathbb{F}_2)^2$ is invariant under $\SL_2(\mathbb{F}_2)$, but one can check by hand that there is no nonzero $\SL_2(K)$-invariant quadratic form on $K^2$ for any infinite field $K$.  So there is no nonzero $\SL_2$-invariant quadratic form on $k^2$ for any field $k$.
Obviously the distinction between $G$ and $G(k)$-invariance exists only for $k$ finite.
\end{borel*}

\begin{borel*} \label{bil.basic}
The usual proofs for irreducible representations in characteristic zero from \cite{Bou:g7} or \cite{GW} show that a Weyl module $V(\la)$ has a nonzero $G$-invariant bilinear form if and only if $\la = -w_0 \la$ for $w_0$ the longest element of the Weyl group, and that a nonzero $G$-invariant bilinear form is uniquely determined up to a factor in $\kx$.  (The uniqueness comes from the fact that $\End_G(V(\la))$ is $k$ \cite[II.2.13]{Jantzen}.)

Fix a surjection $\pi \!: V(\la) \ra L(\la)$.  Its kernel is the \emph{(representation-theoretic) radical} of $V(\la)$, i.e., the intersection of the maximal proper submodules \cite[II.2.14]{Jantzen}.  Every nonzero $G$-invariant bilinear form $\sbar$ on $L(\la)$ gives rise to a nonzero $G$-invariant bilinear form $s$ on $V(\la)$ by setting $s(x,y) := \sbar(\pi(x), \pi(y))$.  As $\sbar$ is necessarily nondegenerate, the radical of $s$ equals the kernel of $\pi$.  By uniqueness, the radical of every nonzero $G$-invariant bilinear form on $V(\la)$ is $\ker \pi$.
\end{borel*}

\begin{borel*}
In characteristic 2, quadratic forms---not symmetric bilinear forms---are the objects naturally associated with algebraic groups of types $B$ and $D$ \cite[\S3]{Ti:C}.  Every quadratic form $q$ gives a symmetric bilinear form $b_q$ by setting $b_q(x, y) := q(x +y) - q(x) - q(y)$.  The \emph{radical} of $q$ is defined to be the radical of $b_q$.
\end{borel*}

\begin{prop} \label{q.char2}
Let $q$ be a nonzero $G$-invariant quadratic form on $V(\la)$.  Then:
\begin{enumerate}
\item If $\la \ne 0$ or $\chr k \ne 2$, the radical of $q$ is the representation-theoretic radical of $V(\la)$.
\item If $q'$ is a $G$-invariant quadratic form on $V(\la)$, then $q' = cq$ for some $c \in k$.
\end{enumerate}
\end{prop}

Note that (1) fails if $\la = 0$ and $k$ has characteristic 2.  In that case, $V(\la)$ is the irreducible 1-dimensional representation, hence the radical of $V(\la)$ is zero.

\begin{proof}[Proof of Prop.~\ref{q.char2}]
In characteristic $\ne 2$, the proposition is trivial, so we assume that $k$ has characteristic 2.  To prove (1), we may enlarge $k$ and so assume that every element of $k$ is a square.  For sake of contradiction, suppose that (1) is false, which by \ref{bil.basic} is equivalent to $b_q$ being identically zero.  Then $q(x+y) = q(x) + q(y)$ for all $x,y \in V(\la)$ and the elements of $V(\la)$ on which $q$ vanishes form a subspace $W$ invariant under $G$.

We claim that $\dim V/W$ is at most 1.  Otherwise, there are elements $v_1, v_2 \in V(\la)$ whose images in $V/W$ are linearly independent; in particular, $q(v_1), q(v_2)$ are not zero.
Multiplying $v_1$ by an element of $\kx$, we may assume that $q(v_1) = q(v_2)$, so that $v_1 + v_2$ belongs to $W$.  This contradicts the linear independence, proving the claim.

If $W$ is a proper subspace, then it is contained in the radical of $V(\la)$, hence maps to zero in $L(\la)$.  It follows that the dimension of $V/W$ is at least the dimension of $L(\la)$, which is at least 2 because $\la \ne 0$.  We conclude that $W$ is all of $V$, i.e., that $q$ is zero, a contradiction.

Claim (2) follows easily from the proof of (1), as uniqueness for bilinear forms on $L(\la)$ gives a $c \in k$ such that $c b_q - b_{q'}$ is the zero form.
\end{proof}

We have treated uniqueness for and the radicals of $G$-invariant quadratic forms, but what of existence?  Consider the following:
\begin{eg} \label{SL2.L}
For a nonnegative integer $n$, \emph{the irreducible representation $L(n)$ of $\SL_2$ has a nonzero $\SL_2$-invariant quadratic form if and only if $n$ is even ($\chr k \ne 2$) or $n$ is not a power of $2$ ($\chr k = 2$).}  To prove the characteristic 2 result for $n \ne 0$, one writes $L(n)$ as a tensor product $\ot L(1)^{[j]}$ where the product runs over those $j$ such that the term $n_j$ in the base-2 expansion of $n$ is not zero and $L(1)^{[j]}$ denotes the $j$-th Frobenius twist of $L(1)$, i.e., the composition of the standard 2-dimensional representation of $\SL_2$ with the map $\SL_2 \ra \SL_2$ induced by the homomorphism $x \mapsto x^{2^j}$ of $k$.  There is no nonzero $\SL_2$-invariant quadratic form on $L(1)$ by \ref{clar}, hence none on $L(1)^{[j]}$ for all $j$.  Conversely, $L(1)^{[j]}$ has a nonzero $\SL_2$-invariant alternating bilinear form for all $j$, so $\ot L(1)^{[j]}$ has a nonzero $\SL_2$-invariant quadratic form if there are at least two terms in the product by, e.g., \cite[p.~67, Exercise 21]{KMRT} or \cite[3.4]{SinWillems}.
\end{eg}

\begin{borel*} \label{singular}
In general, we have: \emph{There is a nonzero $G$-invariant quadratic form on the irrep $L(\la)$ if and only if there is a nonzero $G$-invariant quadratic form on the Weyl module $V(\la)$ that is zero on the radical of $V(\la)$.}  Indeed, you can lift a $G$-invariant quadratic form on $L(\la)$ to one on $V(\la)$ that is zero on the radical, so the if direction is trivial and the only if direction follows from Prop.~\ref{q.char2}.
\end{borel*}

\begin{eg} \label{SL2.V}
For a non-negative integer $n$, \emph{the Weyl module $V(n)$ of $\SL_2$ has a nonzero $\SL_2$-invariant quadratic from if and only if $n$ is even ($\chr k \ne 2$) or $n \ne 1$ ($\chr k = 2$).}  To check this in characteristic 2, one constructs a quadratic form on $V(n)$ for $n$ even by reducing a form defined over $\Z$, see the next section.  The claim for $n$ odd follows from the two preceding paragraphs.
\end{eg}

\begin{rmk} \label{V.odd}
In the rest of this paper, we have ignored the Weyl modules $V(n)$ of groups of type $A_1$ with $n$ odd.  We have done so for two reasons.  First, \emph{these modules are defined only for $\SL_2$ and $\PGL_2$}, i.e., they are defined for $\SL(Q)$ only if $Q$ is split.  Indeed, suppose that there is a representation $W$ of $\SL(Q)$ defined over $k$ that is equivalent to $V(n)$ over a separable closure $\ksep$ of $k$.  The radical of $V(n)$ is a characteristic submodule, so by Galois descent, there is a $k$-defined submodule $R$ of $W$ such that $W/R$ is equivalent to the irrep $L(n)$ over $\ksep$.  This implies that $n[Q] = 0$ in the Brauer group of $k$ \cite[3.5]{Ti:R}, hence that $Q$ is split (because $n$ is odd).

Second, if $k$ has characteristic different from 2, the nonzero $\SL_2$-invariant bilinear form on $V(n)$ is skew-symmetric and its isomorphism class is uniquely determined by the fact that its radical is the radical of $V(n)$.  If $k$ has characteristic 2 and $n \ge 3$, then there is a nonzero $\SL_2$-invariant quadratic form on $V(n)$ that vanishes on the radical of $V(n)$ by the preceding example.  This form is always hyperbolic, because $0$ is not a weight of $L(n)$.
\end{rmk}

\begin{borel*} \label{VvL}
The material in this section implies that---for answering the question ``What is $q$?" posed at the beginning of this article for absolutely irreducible representations---it suffices to consider the Weyl modules.  Suppose that you know a 
 $G$-invariant quadratic form $q$ on a Weyl module $V(\la)$.
If the field $k$ has characteristic $\ne 2$, then obviously the invariant form on the irrep $L(\la)$ is the nondegenerate part of $q$.  If $k$ has characteristic 2, one applies the criterion \ref{singular} to determine if $L(\la)$ has a nonzero $G$-invariant quadratic form; if it does, then again it is the nondegenerate part of $q$.
\end{borel*}

\section{Theorems \ref{REP} and \ref{QUAD} in characteristic 2} \label{char2}

In this section, the field $k$ has characteristic 2.  We prove analogues of Theorems \ref{REP} and \ref{QUAD} in this setting.

We adopt the notation of \cite{Baeza} for symmetric bilinear forms and quadratic forms over $k$.  In particular, $\qform{\alpha_1, \ldots, \alpha_n}$ denotes the bilinear form $(x, y) \mapsto x^t A y$ on $k^n$ where $A$ is the diagonal matrix whose $(i,i)$-entry is $\alpha_i$.

Write $s(n)$ for the number of 1's appearing in the binary expansion of the natural number $n$.  We have:
\begin{specB}[$\chr k = 2$]
For $n$ even, the bilinear form $\podd$ is identically zero and the bilinear form $\peven$ has nondegenerate part $2^{s(n) - 1}\qform{1}$.
\end{specB}

The theorem says the nondegenerate part of $\peven$ is $\qform{1}$ if $n$ is a power of 2 and is metabolic otherwise.

\begin{proof}
The binomial coefficient  $\binom{n}{m}$ is even if and only if some 1 in the binary expansion of $m$ appears where there is a 0 in the binary expansion of $n$.  Hence $\binom{n}{m}$ is even for all odd $m$, which proves that $\podd$ is identically zero.  Further, the coefficient $\binom{n}{n/2}$ is even and the number of $\binom{n}{m}$ that are odd is $2^{s(n)}$.
The $2^{s(n)}$ odd binomial coefficients are divided evenly between those even $m$ with $0 \le m < n/2$ and those even $m$ with $n/2 < m \le n$.  Hence the first collection has exactly $2^{s(n) - 1}$ members.  This proves the claim about $\peven$.
\end{proof}

We fix notations in order to state Theorem \ref{REP}.
A quaternion algebra $Q$ over $k$ is generated by elements $i, j$ such that 
\[
i^2 + i = a, \quad j^2 = b, \quad \text{and} \quad ij = ji + j
\]
for some $a \in k$ and $b \in \kx$.  The element $i$ generates a quadratic \'etale $k$-algebra $K$ in $Q$ isomorphic to $k[\alpha] / (\alpha^2 + \alpha - a)$, and the norm from this extension is a 2-dimensional quadratic form denoted by $[1, a]$.  We write $[c]$ for the 1-dimensional space with basis $e$ endowed with the quadratic form sending $e \mapsto c$.  

The norm of $Q$ restricted to the trace zero elements is
\[
Q' := [1] \oplus \qform{b} [1, a].
\]
This can be calculated by identifying $Q \ot K$ with $M_2(K)$ via 
\begin{equation} \label{ch2.Q}
i \mapsto \stbtmat{\alpha}{0}{0}{\alpha+1} \quad \text{and} \quad j \mapsto \stbtmat{0}{b}{1}{0}
\end{equation}
and restricting the determinant of $M_2(K)$ to the $k$-subspace spanned by $1, j, ij$.

Recall that every quadratic form can be written uniquely as $q_0 \oplus q_\nd$, where $q_0$ is the zero quadratic form and $q_\nd$ is nondefective, see e.g.~\cite[2.4]{HL}; we call $q_\nd$ the \emph{nondefective part} of $q$.  We write $\H$ for the hyperbolic plane $[0,0]$.

\begin{specA}[$\chr k = 2$]
Let $n$ be even.  There is an $\SL(Q)$-invariant quadratic form on the Weyl module $V(n)$ that is isomorphic to 
\[
\begin{cases}
Q' & \text{if $n = 2$, and} \\
[1] \oplus \H \oplus (n-2) [0] &\text{if $n$ is a power of $2$ and $n \ge 4$.} \\
\end{cases}
\]
Otherwise, the nondefective part of the $\SL(Q)$-invariant form is $2^{s(n)-1}$ hyperbolic planes.
\end{specA}

The form given by the theorem is uniquely determined (up to a scalar factor) by the property of being $\SL(Q)$-invariant, see Prop.~\ref{q.char2}(2).

\begin{proof}
The quadratic form $v \mapsto f(v,v)$ over $\Z$ defined by \eqref{f.def} is divisible by 2 because $\binom{n}{n/2}$ is even.  Dividing by 2 and extending scalars to $k$, we find that the quadratic form 
\begin{equation} \label{ch2.q}
\left( \frac12 \binom{n}{n/2} \right) x_{n/2}^2 + \sum_{\stackrel{0 \le m < n/2}{\text{$m$ even}}} \binom{n}{m} x_m x_{n-m}
\end{equation}
is $\SL_2$-invariant, where $x_0, x_1, \ldots, x_n$ is a basis of $V(n)^*$ dual to $e_0, e_1, \ldots, e_n$.  (The terms $\binom{n}{m} x_m x_{n-m}$ with $m$ odd are zero as in the proof of Th.~\ref{QUAD}.)  The isomorphism class of \eqref{ch2.q} is 
\[
\left[ \frac12\binom{n}{n/2} \right] \oplus\ (\peven \cdot \H).
\]
The integer $\frac12 \binom{n}{n/2}$ is odd if and only if $n$ is a power of 2 by Kummer's criterion from \S\ref{binomial} or \cite[10.29]{KMRT}.  Applying Th.~\ref{QUAD} completes the proof of Th.~\ref{REP} in case $Q$ is split.

Suppose now that $Q$ is not split, and in particular that the subalgebra $K$ of $Q$ generated by the element $i$ is a separable field extension of $k$.  Write $\iota$ for the non-identity $k$-automorphism of $K$, i.e., the map such that $\iota(\alpha) = \alpha + 1$.  Via \eqref{ch2.Q}, $Q$ is identified with the $k$-subalgebra of $M_2(K)$ fixed by $\Int(\eta_\iota) \iota$ for $\eta_\iota$ as in \eqref{eta.def}.  The same computations as in \S\ref{desc} show that the $\SL(Q)$-invariant form on $V(n)$ has the same nondefective part as
\begin{equation} \label{ch2.q1}
\left[ \frac12 \binom{n}{n/2} \right] \oplus \peven \cdot \qform{b^{n/2}} \cdot [1, a].
\end{equation}
This verifies the theorem in case $n = 2$.  If $n$ is not a power of 2, we are done by Th.~\ref{QUAD} because `metabolic tensor quadratic is hyperbolic' \cite[I.4.7]{Baeza}.  

Finally, if $n$ is a power of 2 and not 2, the nondefective part of \eqref{ch2.q1} is $[1] \oplus [1, a]$.  But this is the norm on the trace zero elements of the quaternion algebra with $b = 1$, i.e., on $M_2(k)$, hence this form is isomorphic to $[1] \oplus \H$.
\end{proof}

\medskip
\noindent{\small{\textbf{Acknowledgment.} I thank George McNinch and Pham Tiep for their helpful remarks on representations and Jean-Pierre Serre for suggesting the problem and for his comments on an earlier version of this note.}}

\providecommand{\bysame}{\leavevmode\hbox to3em{\hrulefill}\thinspace}
\providecommand{\MR}{\relax\ifhmode\unskip\space\fi MR }
\providecommand{\MRhref}[2]{%
  \href{http://www.ams.org/mathscinet-getitem?mr=#1}{#2}
}
\providecommand{\href}[2]{#2}

\end{document}